\documentstyle[12pt]{article}

\newlength{\abstractwidth}
\renewcommand{\thefootnote}{\fnsymbol{footnote}}
\renewcommand{\thanks}[1]{\footnote{#1}} 
\newcommand{\starttext}{ \setcounter{footnote}{0}
\renewcommand{\thefootnote}{\arabic{footnote}}}

\newcommand{\be}{\begin{equation}}
\newcommand{\bea}{\begin{eqnarray}}
\newcommand{\eea}{\end{eqnarray}} 
\newcommand{\ee}{\end{equation}}

\def\ba{\begin{eqnarray}}
\def\ea{\end{eqnarray}}

\def\o{\omega}

\def\det{{\rm det}}

\def\log{\,{\rm log}\,}

\def\o{\omega}

\def\o{\omega}

\def\na{\nabla}
\def\p{\partial}

\def\C{{\bf C}}

\def\ddb{{\partial\bar\partial}}

\def\na{{\nabla}}

\def\[{{\bf [}}
\def\]{{\bf ]}}

\begin{document}
\starttext \baselineskip=15pt \setcounter{footnote}{0}
\newtheorem{theorem}{Theorem}
\newtheorem{lemma}{Lemma}
\newtheorem{definition}{Definition}
\newtheorem{proposition}{Proposition}
\newtheorem{corollary}{Corollary}

\begin{center}
{\Large \bf
NEW CURVATURE FLOWS IN COMPLEX GEOMETRY
 \footnote{Contribution to the proceedings of the Conference celebrating ``50 Years of the Journal of Differential Geometry", Harvard University, April 28-May 2, 2017. Work supported in part by the National Science Foundation Grants DMS-12-66033 and DMS-1605968, and the Simons Collaboration Grant-523313.}}
\end{center}

\centerline{Duong H. Phong, Sebastien Picard, and Xiangwen Zhang}

\bigskip

\begin{abstract}

{\small The Anomaly flow is a flow of $(2,2)$-forms on a $3$-fold which was originally motivated by string theory and the need to preserve the conformally balanced property of a Hermitian metric in the absence of a $\p\bar\p$-Lemma.
It has revealed itself since to be a remarkable higher order extension of the Ricci flow. It has also led to several other curvature flows which may be interesting from the point of view of both non-K\"ahler geometry and the theory of non-linear partial differential equations. This is a survey of what is known at the present time about these new curvature flows.}

\end{abstract}

\section{Introduction}
\setcounter{equation}{0}

The search for canonical metrics is at the interface of several particularly active fields in mathematics. The theory of non-linear partial differential equations certainly plays a major role, since a canonical metric is typically defined by a curvature condition, and hence is a non-linear system of second order equations in the metric. The curvature condition should optimize the metric in some sense, which suggests a variational principle, or a link with the field equation of a physics theory. The existence of a canonical metric, or of a minimizer for the variational problem, is usually tied to some deep geometric properties of the underlying space, such as stability in algebraic geometry. This interaction between the theory of partial differential equations, theoretical physics, and algebraic geometry has been particularly fertile for complex, and especially K\"ahler geometry. There, the complex and the Riemannian structures fit seamlessly, and the search for Hermitian-Yang-Mills connections \cite{D,UY}, K\"ahler-Einstein metrics \cite{Aubin, Y, CDS3}, and K\"ahler metrics with constant scalar curvature (see \cite{D1, CC} for recent advances and \cite{PS} for a survey) have motivated many exciting developments in the last 50 years.

\medskip
However, conditions such as Yang-Mills connections, or metrics of constant scalar or Ricci curvature are all linear in the curvature form. Recent
advances in theoretical physics, especially string theory, suggest that the notion of canonical metrics should be widened to include conditions which involve its higher powers. This is in view of the Green-Schwarz anomaly cancellation mechanism
\cite{GS}, which involves the square of the curvature, and also of string-theoretic corrections to the Einstein-Hilbert and the Yang-Mills actions. Since the curvature form with respect to the Chern connection of the Hermitian metric is a $(1,1)$-form (valued in the space of endomorphisms of the tangent space or a vector bundle), equations involving its $p$-th power for $p>1$
will be equations for $(p,p)$-forms. But the K\"ahler condition is a condition on $(1,1)$-forms, so we expect to have to replace it by conditions on $(p,p)$-forms, such as Michelsohn's balanced condition \cite{Michelsohn}, which will usually be weaker. Because of these accumulated difficulties, systems involving powers of the curvature may seem inaccessible but for the breakthrough by Fu and Yau \cite{FY1, FY2}, who found some special solutions by PDE methods to one such system, namely the Hull-Strominger system \cite{Hull1,Hull2,S}. The Fu-Yau solution also revealed that the system has a lot of structure and may perhaps be amenable to a more systematic approach.

\smallskip
If we have to implement a condition which is weaker than the K\"ahler condition, we face a very serious obstacle, namely the lack of a $\p\bar\p$-Lemma in non-K\"ahler geometry. Perhaps surprisingly, in some cases such as the Hull-Strominger system, this can be circumvented by introducing a suitable geometric flow, which is a flow of $(2,2)$-forms on a $3$-dimensional complex manifold \cite{PPZ2}. Even more surprising is that the resulting flow, called the Anomaly flow, turns out to be a generalization of the Ricci flow, albeit with corrections due to the metric not being K\"ahler and to the higher powers of the curvature tensor \cite{PPZ5}. Besides its original motivation as a parabolic approach to solving the Hull-Strominger system, the Anomaly flow admits many variants and generalizations, which can potentially be useful in addressing some classic questions of non-K\"ahler geometry. While many extensions of the Ricci and related flows have been introduced over the years in the literature \cite{BV, Bry, BryXu, Gill, Ka, LYZ, LoWe, ST, ST2, TW2}, the ones arising from the Anomaly flow seem all new. In fact, while for many of the existing generalizations, only short-time existence theorems and/or derivative estimates are available at the present time,
the long-time behavior and convergence of the Anomaly-related flows can be established in a number in highly non-trivial cases. This can perhaps be attributed to their natural handling of the conformally balanced condition, and also argues for their interest from the point of view of the theory of non-linear partial differential equations.

\medskip
The purpose of the present paper is to present an introduction to the Anomaly and related flows, and to survey some of what is known.
Their study has just begun, and many basic questions are as yet unanswered. We hope that the paper will motivate many researchers to help answer these questions.

\section{The Reference Flow: the K\"ahler-Ricci flow}
\setcounter{equation}{0}

When discussing the new flows, it will be instructive to compare them to a reference flow in complex geometry, namely the K\"ahler-Ricci flow. We begin by recalling some basic facts about the K\"ahler-Ricci flow. The literature on the subject is immense, and we shall only touch upon those aspects where a comparison with the new flows would be particularly helpful.

\subsection{The Ricci flow}
We start with the Ricci flow. Let $X$ be a Riemannian manifold. For simplicity, we shall assume throughout this paper that $X$ is compact, unless indicated specifically otherwise. The Ricci flow is the flow of metrics $t\to g_{ij}(t)$ given by
\bea
\p_t g_{ij}(t)=-2\,R_{ij}(t), \qquad g_{ij}(0)=g^0_{ij}
\eea
where $g^0_{ij}$ is an initial metric, and $R_{ij}(t)$ is the Ricci curvature of $g_{ij}(t)$. In normal coordinates near any given point, the Ricci curvature is given by
\bea\label{RF-heat}
R_{ij}= {1 \over 2} ( -g^{pq}\p_p\p_q\,g_{ij}- g^{pq} \p_i\p_j \, g_{pq} + g^{pq} \p_p\p_j\, g_{iq} + g^{pq} \p_i\p_q \, g_{pj}).
\eea
Using De Turck's trick \cite{De}, the last three terms can be removed by gauge fixing, and the first term $-g^{pq}\p_p\p_q\,g_{ij}$ is used to obtain a unique short time solution to the Ricci flow.

\smallskip
Initially developed by Hamilton, the early successes of the Ricci flow include the proof of its convergence to a space form for $3$-folds with positive Ricci curvature 
\cite{Hamilton1} and, with appropriate surgeries, for $4$-folds with positive curvature operator \cite{Hamilton2}. As everyone knows, these successes culminated in Perelman's spectacular proof of the Geometrization Conjecture \cite{Perelman1,Perelman2}. It has been a very active research direction ever since, and many deep and powerful results continue to be obtained.

\subsection{The K\"ahler-Ricci flow}

Let $X$ be now a compact K\"ahler manifold. The K\"ahler-Ricci flow is just the Ricci flow
\bea
\partial_t g_{\bar kj}=-R_{\bar kj}(t)
\eea
where $g_{\bar kj}(0)$ (and hence $g_{\bar kj}(t)$ for all $t$) is K\"ahler.

\medskip

We recall the notion of Chern unitary connection defined by a Hermitian metric $g_{\bar kj}$. To a Hermitian metric $g_{\bar kj}$, we can associate its form $\omega=ig_{\bar kj}dz^j\wedge d\bar z^k$. The metric is said to be K\"ahler if $d\omega=0$. The K\"ahler class of $\omega$ is then the de Rham cohomology class $[\omega]$.

Given a Hermitian metric $\omega$, the Chern unitary connection $\na$ is defined by the requirement that $\nabla_{\bar k}V^j=\partial_{\bar k}V^j$ and unitarity, i.e., $\nabla_{\bar k}g_{\bar mj}=0$.
In general, the Chern unitary connection differs from the Levi-Civita connection. But the two are the same if $\omega$ is K\"ahler.

\smallskip
For K\"ahler metrics, the Ricci curvature is given by the very simple formula
\bea
\label{ricci}
R_{\bar kj}=-\partial_j\partial_{\bar k}\log\,\o^n,
\eea
where $n$ is the complex dimension of $X$. The corresponding Ricci form $Ric(\o)$ is defined by $Ric(\omega)=iR_{\bar kj}dz^j\wedge d\bar z^k$, or equivalently,
\bea
Ric(\omega)=-i\partial\bar\partial \log \omega^n.
\eea
Clearly $Ric(\omega)$ is a closed form, and defines itself a cohomology class.
Since for two different metrics $\omega$ and $\tilde\omega$, the ratio $\omega^n/\tilde\omega^n$ is a well-defined smooth function, $[Ric(\omega)]$ is independent of $\omega$. It is called the first Chern class $c_1(X)$. More geometrically, let $K_X$ be the canonical bundle of $X$, that is, the bundle of $(n,0)$-forms on $X$. We can think of $\o^n$ as a positive section of $K_X\otimes\bar K_X$, and hence as a metric on the anti-canonical bundle $K_X^{-1}$. The Ricci curvature defined by (\ref{ricci}) can then be viewed as the curvature of the line bundle $K_X^{-1}$, and its de Rham cohomology class $c_1(X)$ can be identified as the first Chern class of $K_X^{-1}$.

\subsection{Immediate advantages of the K\"ahler condition}

Since $\p_t\o=-Ric(\o)$, the closedness of $Ric(\o)$ implies that $\o$ evolves only by closed forms, and hence it remains closed if it is closed at time $t=0$. In other words, the K\"ahler property is preserved. The resulting cohomological information accounts for some marked differences between the K\"ahler-Ricci flow and the general Ricci flow on just Riemannian manifolds.

\smallskip

(a) In the K\"ahler case, the K\"ahler-Ricci flow implies
$\partial_t[\omega(t)]=-c_1(X)$, and hence, even though $\o(t)$ may still be difficult to find, we know at least its cohomology class
\bea
[\omega(t)]=[\omega(0)]-t\,c_1(X).
\eea

(b) Thus the flow can only exist for $t$ such that $[\omega(t)]>0$, in the sense that $[\omega(t)]$ contains a positive representative. This suggests that the maximum time of existence for the flow should be given by
\bea
T={\rm sup}\{t\,:\, [\omega(0)]-t\,c_1(X)>0\}
\eea
which is determined by the cohomology class of $[\omega(0)]$, and does not depend on $\omega(0)$ itself. This was confirmed  by Tian and Z. Zhang \cite{TZ}.

\smallskip

(c) On a compact K\"ahler manifold, the $\partial\bar\partial$-Lemma says that if $\o$ and $\tilde\o$ are any two $(1,1)$-forms satisfying $[\tilde\omega]= [\omega]$, then we can write $\tilde\omega=\omega+i\partial\bar\partial\varphi$ for some smooth function $\varphi$ unique up to an additive constant, which is called the potential of $\tilde\omega$ with respect to $\o$. In particular, a K\"ahler metric $\o$ is determined by its cohomology class $[\o]$ and a scalar function. We can apply the $\p\bar\p$-Lemma to reduce the K\"ahler-Ricci flow to a scalar flow as follows. First we choose a known and evolving $(1,1)$-form $\hat\o(t)$ with $\p_t\hat\o(t)=i\p\bar\p \log\Omega=-c_1(X)$, $\hat\o(0)=\o(0)$, where $\Omega$ is a fixed smooth, strictly positive section of $K_X$. Then $[\hat\o(t)]=[\o(t)]$, and we can set
\bea
\o(t)=\hat\o(t)+i\p\bar\p\varphi(t)
\eea
with $\varphi(t)$ now the unknown scalar function. The K\"ahler-Ricci flow is equivalent to
\bea
i\p\bar\p\p_t\varphi(t)=i\p\bar\p\log \,\o^n(t)-i\p\bar\p\log\Omega
=i\p\bar\p \log{\o^n(t)\over\Omega}
\eea
Changing if necessary $\Omega$ by a constant multiple, this equation is equivalent to the following parabolic complex Monge-Amp\`ere type equation for the scalar function
$\varphi$,
\bea
\p_t\varphi=\log{(\hat\o(t)+i\p\bar\p \varphi)^n\over\Omega},
\qquad
\hat\o(t)+i\p\bar\p\varphi>0.
\eea
This equivalence between a flow of metrics and a much simpler flow of scalar functions is arguably the most useful consequence of the K\"ahler property, and one which will be sorely missing when we consider the new curvature flows.

\subsection{The K\"ahler-Ricci flow in analytic/algebraic geometry}

Early on, the K\"ahler-Ricci flow provided a parabolic approach to the problem of finding K\"ahler-Einstein metrics, that is K\"ahler metrics satisfying the condition \bea
Ric(\omega)=\mu\,\omega
\eea
for $\mu$ constant. This problem was solved using elliptic methods by Yau \cite{Y} for $\mu=0$ and by Yau \cite{Y} and Aubin \cite{Aubin} for $\mu<0$. Such metrics can also be viewed as stationary points of the (normalized) K\"ahler-Ricci flow, and an alternative proof using the K\"ahler-Ricci flow was given by H.D. Cao \cite{Cao}. 
The case $\mu>0$ was solved only relatively recently in 2013 by Chen, Donaldson, and Sun \cite{CDS1, CDS2, CDS3}. Subsequent progress on the K\"ahler-Ricci flow approach in this case can be found in \cite{TianZhang, CW} and references therein.

\smallskip

In the previous application, the initial assumption, which is necessary for the existence of K\"ahler-Einstein metrics, was that $c_1(X)$ was either zero, or definite. More recently, in an important development, the existence of K\"ahler-Einstein metrics was also established by Wu and Yau \cite{WY1,WY2, WY3}, starting instead from an assumption of negative holomorphic sectional curvature. 
Wu and Yau used elliptic methods, but  Nomura \cite{Nomura}, Tong \cite{Tong}, and Huang et al \cite{HLTT} have been able to recover and extend the results of Wu and Yau in several ways, using the K\"ahler-Ricci flow.

\smallskip 

The above applications of the K\"ahler-Ricci flow concern the cases when a canonical metric does exist, and the K\"ahler-Ricci flow converges to it.
Several longer range and presently ongoing programs are based on the expectation that the K\"ahler-Ricci flow can still help find a canonical metric, or a canonical complex structure, even when there is none on the original manifold. These include 
the analytic minimal model program of Song and Tian \cite{SongTian}, with early contributions by Tsuji \cite{Ts}, 
which requires the continuation of the K\"ahler-Ricci flow beyond singularities, possibly on a different manifold; and flows on instable manifolds,
where one would have to address the phenomena of jumps to a different complex structure and moduli theory. A major difficulty in general is to identify the limiting complex manifold. So far only the case of $S^2$ with an arbitrary configuration of marked points, including semistable and unstable ones, has been fully worked out in \cite{PSSW}. The convergence in the stable case had been obtained in \cite{MRS}.

\section{Motivation for the New Flows}
\setcounter{equation}{0}

While the study of the K\"ahler-Ricci flow remains a very active research area, with many hard and deep questions left open, a number of new problems arising from different areas seem to require new flows. Remarkably, while their origin is a priori unrelated, these new flows turn out to be in many ways
just higher order corrections of the Ricci flow, in a natural setting which is somewhere intermediate between the Riemannian and the K\"ahler settings. We begin with an informal discussion of the motivation and qualitative features of these flows, leaving the precise formulation to the next section.

\subsection{Anomaly cancellation and powers of the curvature}

\smallskip

In string theory, a crucial requirement discovered by M. Green and J. Schwarz \cite{GS} for the cancellation of anomalies is the equation
\bea
\label{anomaly}
dH
={\alpha'\over 4}\left({\rm Tr}(Rm \wedge Rm)-{\rm Tr}(F\wedge F)\right).
\eea
Here space-time is a $10$-dimensional manifold, and the fields in the theory include a metric $g_{ij}$, a Hermitian metric $h$ on a vector bundle $E$, 
and a $3$-form $H$. The expressions
$Rm$ and $F$ are the curvatures of $\omega$ and $h$, viewed as $2$-forms valued in the endomorphisms of the tangent space and the bundle $E$ respectively. The coefficient $\alpha'$ is the string tension, which is a positive constant.

\smallskip
Note that this is an equation of $4$-forms involving the square of the curvature. By contrast, more familiar equations such as the K\"ahler-Einstein equation or the Yang-Mills equation are linear in the curvature. The right hand side of the equation (\ref{anomaly}) arises rather from Pontryagin or Bott-Chern characteristic classes, which arise themselves from gauge and gravitational anomalies. To derive it from an action, one would need to expand in the string tension $\alpha'$, and keep both the leading terms as well as the next order expansion terms.

\smallskip

As shown by Candelas, Horowitz, Strominger, and Witten \cite{CHSW}, in compactifications of the heterotic string from $10$-dimensional to $4$-dimensional space-time, supersymmetry requires the existence of a covariantly constant spinor from which an integrable complex structure with nowhere vanishing holomorphic $(3,0)$-form $\Omega$ can be constructed on the internal $6$-dimensional manifold. Restricted to the internal space, the metric $g_{ij}$ can then be identified with a $(1,1)$-form $\o$, and $H$ becomes the torsion of $\o$, $H=i(\bar\p-\p)\o/2$. The bundle $E$ turns out to be a holomorphic vector bundle, and the curvature $F$ is that of the Chern unitary connection of $h$ on $E$. In the equation (\ref{anomaly}), it is natural to view $Rm$ as the curvature of the Chern unitary connection of $\o$. However, it has been shown by Hull \cite{Hull2} that other unitary connections are also admissible. Altogether, the equation (\ref{anomaly}) is then an equation for $(2,2)$-forms, which implies that ${\rm Tr}(Rm\wedge Rm)$ is a $(2,2)$-form, although $Rm$ itself may not be a $(1,1)$-form.

\smallskip
We observe that, except in the special case when the right hand side of (\ref{anomaly}) vanishes identically (which can only happen if the vector bundles $T^{1,0}(X)$ and $E$ have the same second Chern class), the equation (\ref{anomaly}) implies that the Hermitian metric $\o$ is not K\"ahler. Thus, besides introducing an equation of $(2,2)$-forms and the square of the curvature tensor, the anomaly cancellation mechanism of Green and Schwarz also suggests a widening of metrics under consideration to non-K\"ahler metrics.

\subsection{Weakenings of the K\"ahler condition}

Many weakenings of the K\"ahler condition on Hermitian metrics are known, and we shall actually discuss some of them in section \S 6 below. But for our purposes, the most important condition is the notion of balanced metric introduced by Michelsohn \cite{Michelsohn}. A metric $\o$ on an $n$-dimensional complex manifold 
is said to be {\it balanced} if
\bea
\label{balanced}
d(\omega^{n-1})=0
\eea
and a manifold is said to be balanced if it admits a balanced metric. Depending on the circumstances, this condition may well prove to be
preferable  to the K\"ahler condition, because it is invariant under birational transformations, as was shown by Alessandrini and Bassanelli \cite{AB}.

\smallskip

Perhaps unexpectedly, the notion of balanced metric arose relatively recently in a
completely different context, which is the one discussed in the previous section of compactifications of the heterotic string to $4$-dimensional Minkowski space preserving supersymmetry. In the original set-up considered in \cite{CHSW}, the compactification was achieved by a product of Minkowski space time with the internal space. In the generalization considered by
Hull \cite{Hull1, Hull2} and Strominger \cite{S}, the product is replaced by a warped product. Supersymmetry still leads to the internal space being a Calabi-Yau 3-fold, with a holomorphic non-vanishing (3, 0)-form $\Omega$, but the Hermitian metric $\omega$ must satisfy the condition
\bea
\label{c-balanced}
d(\|\Omega\|_\omega\omega^2)=0
\eea
(as reformulated by Li and Yau \cite{LY}). This condition just means that the metric $\|\Omega\|_\o^{1\over 2}\o$ is balanced, and we shall just refer to $\o$ as being {\it conformally balanced}. For more recent developments on balanced metrics, please see the survey paper \cite{Fu} and references therein.

\subsection{Detection of instability}
\setcounter{equation}{0}

\smallskip

Even when the solution of an elliptic equation can already be found by a particular flow under e.g. a stability condition, other flows can be potentially useful too. 
This is because, in the absence of a stability condition, they may fail to converge in different ways, which would provide then different ways of detecting instability. A well-known example is the method of continuity for the K\"ahler-Einstein problem for Fano manifolds. The failure of convergence is detected there by the emergence of a non-trivial multiplier ideal sheaf \cite{Nadel, Siu}. By contrast, the behavior of the K\"ahler-Ricci flow in the unstable case is still obscure. On the other hand, a third flow again with K\"ahler-Einstein metrics as stationary point is the inverse Monge-Amp\`ere flow recently introduced by Collins et al \cite{Co1}. 
This flow appears to provide more information than the K\"ahler-Ricci flow when the underlying manifold is an unstable toric Fano manifold, as it produces then an optimal destabilizing test configuration.

\subsection{New Types of Partial Differential Equations}

\smallskip
The appearance of powers of the curvature, in contrast with the equations of constant scalar or constant Ricci curvature with which we are more familiar,
can be expected to result also in new types of PDE's that have not been encountered before. 

\smallskip
Two new difficulties will also have to be addressed:

\smallskip

(a) The first is the lack of a $\partial\bar\partial$-lemma in non-K\"ahler geometry. So the equations will usually not reduce to scalar equations. In this respect, the non-K\"ahler flows that we shall encounter may be closer to the Ricci flow than the K\"ahler-Ricci flow.

\smallskip

(b) But even when they happen to be reducible to a scalar equation, more often than not, these equations will not be concave in the eigenvalues of the Hessian. This will prevent the use of classical and powerful PDE tools such as those of Caffarelli-Nirenberg-Spruck \cite{CNS, CKNS} and the Evans-Krylov theorem, and other methods will have to be developed.

\smallskip
Thus the new geometric flows may be interesting in their own right from the point of view of the theory of non-linear partial differential equations.

\section{The Anomaly Flow}
\setcounter{equation}{0}

We now formulate precisely the first of the new flows which appeared recently in complex geometry. This is the Anomaly flow introduced by the authors in \cite{PPZ2}, and given this name in recognition of the anomaly cancellation mechanism introduced in \cite{GS}. There are actually several versions of the Anomaly flow, and we begin with the simplest version.

\smallskip

Let $X$ be a compact complex $3$-fold, equipped with a nowhere vanishing holomorphic $(3,0)$-form $\Omega$. Let $t\to \Phi(t)$ be a given path of
closed $(2,2)$-forms, with $[\Phi(t)]=[\Phi_2(0)]$ for each $t$. Let $\omega_0$ be an initial metric which is conformally balanced, i.e., 
$\|\Omega\|_{\omega_0}\omega_0^2$ is a closed $(2,2)$-form. Then the Anomaly flow is the flow of $(2,2)$-forms defined by
\bea
\label{af1}
\partial_t(\|\Omega\|_\omega\omega^2)
=
i\partial\bar\partial\omega-{\alpha'\over 4}
({\rm Tr}\,(Rm\wedge Rm)-\Phi),
\eea
where $\|\Omega\|=(i\Omega\wedge\bar\Omega\,\omega^{-3})^{1\over 2}$ is the norm of $\Omega$ with respect to $\omega$. Here $\alpha'$ is the string tension introduced earlier in (\ref{anomaly}). Mathematically, the equation (\ref{af1}) is well-defined for any real number $\alpha'$ and we fix $\alpha'\in {\bf R}$ unless specified otherwise.

\smallskip

The expression $Rm$ is the curvature of the Chern unitary connection defined by the Hermitian metric $\omega$, viewed as a section of $\Lambda^{1,1}\otimes End(T^{1,0}(X))$. We shall on occasion consider other connections too, and say explicitly so when this is the case.

\smallskip
Note that the stationary points of the flow are precisely given by the Green-Schwarz anomaly cancellation equation (\ref{anomaly}), and 
that the flow preserves the conformally balanced condition. Indeed, by Chern-Weil theory, its right hand side is a closed $(2,2)$-form and hence
for all $t$,
$$
d(\|\Omega\|_\omega\omega^2)=0
$$
at all times $t$, if $d(\|\Omega\|_\omega\omega^2)=0$ at time $t=0$.

\subsection{The Hull-Strominger system}

\smallskip

The main motivation for the Anomaly flow is to solve the following system of equations for supersymmetric compactifications of the heterotic string, proposed independently by Hull \cite{Hull1, Hull2} and Strominger \cite{S}, which generalizes an earlier proposal by Candelas, Horowitz, Strominger, and Witten \cite{CHSW}.
Let $X$ be a $3$-fold equipped with a holomorphic non-vanishing $(3,0)$-form $\Omega$ and a holomorphic vector bundle $E\to X$ with $c_1(E)=0$. Then we look for a Hermitian metric $\omega$ on $X$ and a Hermitian metric $h$ on $E$ on $X$ satisfying
\bea
\label{Hull-Strominger}
i\partial\bar\partial \omega-{\alpha'\over 4}
\left({\rm Tr}(Rm \wedge Rm)-{\rm Tr}(F\wedge F)\right)&=&0\nonumber\\
d(\|\Omega\|_\omega\omega^2)&=&0\nonumber\\
\omega^2\wedge F&=&0
\eea
The third condition is the familiar Hermitian-Yang-Mills equation, so the novelty of the Hull-Strominger system resides essentially in the first two equations.
Actually, the second equation in the above system was originally written in a different way in \cite{Hull1, Hull2} and \cite{S}. The formulation given above, which brings to light the conformally balanced condition, is due to Li and Yau \cite{LY}. 

\smallskip
The solutions of the Hull-Strominger system can be viewed as generalizations of Ricci-flat K\"ahler metrics in the following sense. Assume that 
\bea
\label{anomaly-cancellation}
{\alpha'\over 4}
\left({\rm Tr}(Rm \wedge Rm)-{\rm Tr}(F\wedge F)\right)=0.
\eea
Then the first equation reduces to $i\p\bar\p\o=0$. Combined with the second equation $d(\|\Omega\|_\omega\omega^2)=0$, it implies that $\o$ is both K\"ahler and Ricci-flat (a detalled proof of this fact is given in \S \ref{balance-atheno-kahler} below). One simple way of guaranteeing the condition (\ref{anomaly-cancellation}), as suggested in \cite{CHSW}, is to take $E=T^{1,0}(X)$, and $\o=h$. If $\o$ is K\"ahler, then the first equation in the Hull-Strominger system is satisfied.
The third condition on $F$ coincides then with the Ricci-flat condition on $\o$.
This implies that $\|\Omega\|_\o$ is constant, and hence the second equation is a consequence of $\o$ being K\"ahler. Thus Ricci-flat K\"ahler metrics provide a consistent solution to the Hull-Strominger system. 

\medskip
An immediate difficulty when trying to solve the Hull-Strominger system in general is how to implement the second condition in the absence of an analogue of a $\p\bar\p$-Lemma. While there are many Ansatze that can produce a balanced metric starting from a given one, none of them is particularly special, and they all lead to very unwieldy expressions for their curvatures. Anomaly flows are designed to circumvent this difficulty, and produce solutions of the Hull-Strominger system without appealing to any particular ansatz. The prototype of an Anomaly flow was introduced in (\ref{af1}). For the purpose of solving the Hull-Strominger system, we 
fix a conformally balanced metric $\omega(0)$ on $X$ and a metric $h(0)$ on $E$ and introduce the following coupled flow $t\to (\omega(t),h(t))$,
\bea
\label{af2}
\partial_t(\|\Omega\|_\omega\omega^2)
&=&i\partial\bar\partial \omega-{\alpha'\over 4}
\left({\rm Tr}(Rm \wedge Rm)-{\rm Tr}(F\wedge F)\right)
\nonumber\\
h(t)^{-1}\partial_th(t)
&=&{\omega^2\wedge F\over \omega^3}
\eea
As noted before, the flow preserves the conformally balanced condition for $\o$, and hence its stationary points are automatically balanced if they exist. Since they will manifestly also satisfy the first and third equations in the Hull-Strominger system, they satisfy the complete Hull-Strominger system. Thus we shall have circumvented the problem of the absence of a $\p\bar\p$-Lemma for conformally balanced metrics.

\smallskip
We observe that, while the flow (\ref{af2}) is a coupled flow for $\o(t)$ and $h(t)$, there are cases when $F$ can be viewed as known. For example, in the case of Calabi-Eckmann-Goldstein-Prokushkin fibrations over a $K3$ surface, the Hermitian-Yang-Mills condition on $F$ does not change as $\o$ involves, and $h$ can be taken as fixed, given by the Hermitian-Yang-Mills metric corresponding to the metric $\o(0)$ on $X$.
Setting
$\Phi={\rm Tr}(F\wedge F)$, we see that the flow (\ref{af2}) reduces to the flow (\ref{af1}).

\subsection{Formulation as a flow of metrics instead of $(2,2)$ forms}

While the formulation of the Anomaly flow as a flow of $(2,2)$-forms is ideally suited to the manifest preservation of the conformally balanced condition of Hermitian metrics, for an analytic study of the flow, we need another, more conventional formulation as a flow of $(1,1)$-forms. This is provided by the following theorem \cite{PPZ2, PPZ5}:

\begin{theorem}
\label{af3-theorem}
Consider the anomaly flow (\ref{af1}) with a conformally balanced initial metric.

{\rm (a)} The flow can also be expressed as
\bea
\label{af3}
\partial_t g_{\bar kj}
=
{1\over 2\|\Omega\|_\omega}
\bigg\{-\tilde R_{\bar kj}+g^{s\bar r}g^{p\bar q}T_{\bar qsj}\bar T_{p\bar r\bar k}
-
{\alpha'\over 4}g^{s\bar r}(R_{[\bar ks}{}^\alpha{}_\beta R_{\bar rj]}{}^\beta{}_\alpha-\Phi_{\bar ks\bar rj})\bigg\}
\eea
Here $\tilde R_{\bar kj}=g^{p\bar q}R_{\bar qp\bar kj}$ is the Chern-Ricci tensor, and $T=i\partial\omega={1\over 2}T_{\bar kjm}dz^m\wedge dz^j\wedge d\bar z^k$ is the torsion tensor. The bracket $[,]$ denote anti-symmetrization in each of the two sets of barred and unbarred indices.

{\rm (b)} The flow is parabolic and exists at least for a short time when $\alpha' Rm$ satisfied the following positivity condition
\be
\label{parabolicity1}
|\xi|^2 |A|^2 + {\alpha' \over 2}
g^{q \bar{b}} g^{a \bar{p}} g^{s\bar r} g^{\alpha \bar{\gamma}}  R_{[\bar rq}{}^\beta{}_\alpha \xi_{\bar{p}} \xi_{s]} A_{\bar{\gamma} \beta} A_{\bar{b} a} > 0
\ee 
for all Hermitian tensors $A$ and vectors $\xi$. In particular, it is parabolic and exists at least for a short time if $|\alpha' Rm(0)|<1/2$. 
\end{theorem}

The above formulas for the Anomaly flow should be viewed in the context of conformally balanced metrics. In general, we define the curvature tensor of the Chern unitary connection of a Hermitian metric $g_{\bar kj}$ by
\bea
[\na_j,\na_{\bar k}]V^p=R_{\bar kj}{}^p{}_q V^q
\eea 
for any vector field $V^p$.  Explicitly,
\bea
R_{\bar kj}{}^p{}_q=-\p_{\bar k}(g^{p\bar m}\p_j g_{\bar mq}).
\eea
When $g_{\bar kj}$ is not K\"ahler, there are several candidates for the notion of Ricci curvature,
\bea
R_{\bar kj}=R_{\bar kj}{}^p{}_p,
\quad
\tilde R_{\bar kj}=R^p{}_{p\bar kj},
\quad
R_{\bar kj}'=R_{\bar k}{}^p{}_{pj},
\quad
R_{\bar kj}''=R^p{}_{j\bar kp}
\eea
with corresponding notions of scalar curvature $R,\tilde R, R'$ and $R''$. We always have $R'=R''$ and $R=\tilde R$. However,
for conformally balanced metrics, we also have the following relations
\bea
R_{\bar kj}'=R_{\bar kj}''={1\over 2}R_{\bar kj}
=\tilde R_{\bar kj}-\na^mT_{\bar kjm}.
\eea

For parabolicity purposes, it is crucial that it is the Chern-Ricci curvature $\tilde R_{\bar kj}$ that appears in the right hand side of (\ref{af3}) rather than the other notions of Ricci curvature. Indeed, 
explicitly, 
$$
\tilde R_{\bar kj}=
-g^{p\bar q}g_{\bar km}\partial_{\bar q}( g^{m\bar \ell}\partial_pg_{\bar \ell j})
=-\Delta g_{\bar kj}+\cdots$$
which is the analogue of the equation (\ref{RF-heat}) for the Ricci flow and shows that the flow is parabolic. Note that the other Ricci tensor
$$
R_{\bar kj}=-\partial_j\partial_{\bar k}\log \o^n=
-g^{p\bar q}\partial_j\partial_{\bar k}g_{\bar qp} + \cdots
$$
has a natural interpretation as the first Chern class for $X$, but would not have been appropriate in the right hand side of (\ref{af3}) for parabolicity purposes.

\subsection{Comparison with the K\"ahler-Ricci flow}\label{comparisonRF}

With the formulation (\ref{af3}) for the Anomaly flow, it is easier now to compare it with the K\"ahler-Ricci flow. Perhaps surprisingly, since their original motivations are rather different, the two flows appear rather similar,  with
the Anomaly flow a higher order version of the K\"ahler-Ricci flow
and the additional complications of $\|\Omega\|_\omega^{-1}$, of torsion, and of quadratic terms in the curvature tensor.

\medskip
(a) In the K\"ahler-Ricci flow, the $(1,1)$-cohomology class is determined
$[\omega(t)]=[\omega(0)]-tc_1(X)$. In the Anomaly flow (\ref{af1}), we have rather
$$
\left[\|\Omega\|_{\omega(t)}\omega^2(t)\right]
=
\left[\|\Omega\|_{\omega(0)}\omega^2(0)\right]-t{\alpha'\over 4}\left(c_2(X)-[\Phi_2(0)]\right)$$

\medskip
(b) However, the $(2,2)$ cohomology class $\left[\|\Omega\|_{\omega}\omega^2\right]$ provides much less information than the $(1,1)$-cohomology class $[\omega]$. For example, the volume is an invariant of an $(1,1)$ cohomology class, but not of an $(2,2)$-cohomology class.

\medskip
(c) As an indirect consequence, the maximum time of the Anomaly flow is not determined by cohomology alone and depends on the initial data. In this respect, the Anomaly flow is closer to the Ricci flow than the K\"ahler-Ricci flow.

\medskip
(d) The diffusion operator for the Anomaly flow is more complicated. It includes the Laplacian as for the K\"ahler-Ricci flow, but also the Riemann curvature, e.g.
\bea
\label{curvature-flow}
\delta R_{\bar kj}{}^\rho{}_\lambda
\to
{1\over 2\|\Omega\|}
(\Delta (\delta R_{\bar kj}{}^\rho{}_\lambda)
+2{\alpha'}
g^{\rho\bar\mu}g^{s\bar r}R_{[\bar r\lambda}{}^\beta{}_\alpha
\nabla_s\nabla_{\bar\mu]}\delta R_{\bar kj}{}^\alpha{}_\beta)
\eea
Another complication is that the torsion stays identically $0$ in the K\"ahler-Ricci flow, but it evolves in the Anomaly flow. Explicitly, its evolution is given by
\bea
\p_tT_{\bar pjq}&=&{1\over 2\|\Omega\|_\o}
\bigg[\Delta T_{\bar pjq}
-\alpha'
g^{s\bar r}(\na_j (R_{[\bar p s}{}^\alpha{}_\beta R_{\bar rq]}{}^\beta{}_\alpha
-\Phi_{\bar p s\bar rq})
+\alpha'
g^{s \bar{r}} \na_q(R_{[\bar p s}{}^\alpha{}_\beta R_{\bar rj]}{}^\beta{}_\alpha
-
\Phi_{\bar ps\bar rj}))\bigg]
\nonumber\\
&&
+{1\over 2\|\Omega\|_\o}(T^m{}_{jq}\Psi_{\bar pm}-T_j\Psi_{\bar pq}+T_q\Psi_{\bar pj}+\na_j (T\bar T)_{\bar pq}-\na_q(T\bar T)_{\bar pj})
\nonumber\\
&&
-{1\over 2\|\Omega\|_\o}
(T^r{}_{q\lambda}R^\lambda{}_{r\bar pj}
-
T^r{}_{j\lambda}R^\lambda{}_{r\bar p q})
\eea
where $\Delta=g^{p\bar q}\na_p\na_{\bar q}$ is the Laplacian, and $\Psi$ is the right hand side of (\ref{af1}).

\subsection{The Anomaly flow for the Hull-Strominger-Ivanov system}

In the original formulation of Strominger \cite{S}
of the Hull-Strominger system (\ref{Hull-Strominger}), it was suggested that $Rm$ be the curvature of the Chern unitary connection of $\o$, as $Rm $ would then be a $(1,1)$-form, which would be consistent with the first equation in (\ref{Hull-Strominger}) being an equation of $(2,2)$-forms. On the other hand, in \cite{Hull3}, Hull argued that actually $Rm$ can be the curvature of any connection,
at the only cost of adding finite counterterms to the effective action.

\smallskip

More recently, in \cite{I}, Ivanov argued that, in order to obtain the field equations of the heterotic string, the Hull-Strominger system ought to be supplemented by the condition that $Rm$ be the curvature of an $SU(3)$ instanton. By definition, an $SU(3)$ instanton with respect to a metric $\o$ is given by another Hermitian metric $\tilde\o$ on $X$ whose Chern unitary connection would satisfy
\bea
\label{I}
\o^2\wedge \tilde F=0.
\eea
Here $\tilde F$ is the curvature of the Chern unitary connection of $\tilde\o$, and the right hand side is $0$ because $c_1(X)=0$. Thus a complete system of equations for a supersymmetric vacuum of the heterotic string would consist of three unknowns $(\o,\tilde\o,h)$ satisfying both (\ref{Hull-Strominger}) and (\ref{I}).
We shall refer to the resulting system for $(\o,\tilde\o,h)$ as the Hull-Strominger-Ivanov system. K\"ahler Ricci-flat metrics can still be viewed as special solutions of the Hull-Strominger-Ivanov system, in the sense that if $\o$ is such a metric and we set $E=T^{1,0}(X)$, $h=\tilde\o=\o$, then we would obtain a solution of this system.

\medskip
In \cite{GFRST}, Garcia-Fernandez et al. introduced an elegant way of reformulating the Hull-Strominger-Ivanov system. For this, let $(X,\Omega)$ be as before, and let ${\cal E}$ be a holomorphic vector bundle over $X$. Let $c$ be also a bilinear form
\bea
c: \Lambda^{1,1}({\cal E})\times \Lambda^{1,1}({\cal E})\to \Lambda^{2,2}(X)
\cap {{\rm Ker}\,d}.
\eea
(for some applications, it may be necessary to impose the stronger condition that in the right hand side, ${{\rm Ker}\,d}$ should be replaced by the smaller space
${\rm Im}\,i\p\bar\p$.)
For each Hermitian metric ${\cal H}$ on ${\cal E}$, let $F_{\cal H}$ be the curvature of the corresponding Chern unitary connection. Then the following system of equations for $(\o,{\cal H})$
\bea
\label{GF}
i\p\bar\p \o-{\alpha'\over 4}c(F_{\cal H},F_{\cal H})&=&0\nonumber\\
d(\|\Omega\|_\o\o^2)&=&0\nonumber\\
F_{\cal H}^{0,2}=F_{\cal H}^{2,0}=0,
\quad \o^2\wedge F_{\cal H}&=&0
\eea
is a generalization of the Hull-Strominger-Ivanov system. Indeed, if we let ${\cal E}=T^{1,0}(X)\otimes E$, ${\cal H}=(\tilde\o,h)$, and
\bea
c((A_1,B_1),(A_2,B_2))={\rm Tr}(A_1\wedge A_2)-{\rm Tr}(B_1\wedge B_2)
\eea
where $(A_i,B_i)$ is the corresponding decomposition of $\Lambda^{1,1}({\cal E})
=\Lambda^{1,1}(T^{1,0}(X))\oplus \Lambda^{1,1}(E)$,
we get back the Hull-Strominger-Ivanov system. 
In their paper \cite{GFRST} Garcia-Fernandez et al. also found a variational principle for the Hull-Strominger-Ivanov system, in other words, a space of $(\o,{\cal H})$ on which a functional $M(\o,{\cal H})$ can de defined, whose critical points are exactly the solutions of the system. We won't reproduce the formula for $M(\o,{\cal H})$ here, but note only that a simple version in the special case with $\alpha'=0$ and no ${\cal H}$ appears below in (\ref{M}).

\medskip
It is straightforward to formulate a version of the Anomaly flow for 
the Hull-Strominger-Ivanov system. In view of the previous reformulation (\ref{GF}) of the system, it suffices to consider the following version $t\to (\o(t),{\cal H}(t))$ of the Anomaly flow,
\bea
\label{AF-HSI}
\p_t(\|\Omega\|_\o\o^2)
&=&
i\p\bar\p\o-{\alpha'\over 4} c(F_{\cal H},F_{\cal H})
\nonumber\\
{\cal H}^{-1}\p_t{\cal H}
&=&
{\o^2\wedge F_{\cal H}\over\o^3}.
\eea
for initial data $(\o(0),{\cal H}(0))$ with $\o(0)$ conformally balanced.
Clearly this system again preserves the conformally balanced property of $\o(t)$, and its stationary points satisfy the Hull-Strominger-Ivanov system (\ref{GF}).
In fact, this system is even simpler than the original Anomaly flow defined in (\ref{af2}), because it has no quadratic term in $Rm(\o)$. There are quadratic terms in $F_{\cal H}$, but these are simpler than $Rm(\o)$.

\smallskip

The parabolicity of the flow (\ref{AF-HSI}), and hence its short-time existence is an immediate consequence of the results of \cite{PPZ2} and \cite{PPZ5}, and more specifically of formulas of the type (\ref{af3}).

\section{Illustrative Cases of the Anomaly Flow}
\setcounter{equation}{0}

The Anomaly flow is still largely unexplored, and the many additional difficulties compared to the K\"ahler-Ricci flow which were described in section \S 4.3, suggest that a lot of work will have to be done before we can achieve a level of understanding anywhere close to our present one for the K\"ahler-Ricci flow. However, many special cases of the original Anomaly flow (\ref{af2}) have now been worked out, each of which sheds a different light on the flow, and each of which may ultimately be useful in developing the general theory.

\subsection{The flow on toric fibrations over Calabi-Yau surfaces}

\smallskip

Let $(Y, \hat\omega)$ be a Calabi-Yau surface, equipped with a nowhere vanishing holomorphic (2, 0)-form $\Omega_Y$. Let $\omega_1,\omega_2\in H^2(Y,{\bf Z})$ satisfy $\omega_1\wedge\hat\omega=\omega_2\wedge\hat\omega=0$.

\medskip
From this data, building on ideas going back to Calabi and Eckmann \cite{CE}, Goldstein and Prokushkin \cite{GP} showed how to construct 
a toric fibration $\pi:X\to Y$, equipped with a $(1,0)$-form $\theta$ on $X$ satisfying
$\partial\theta=0$, $\bar\partial\theta=\pi^*(\omega_1+i\omega_2)$. Furthermore,
the form 
$\Omega=\sqrt 3\,\Omega_Y\wedge\theta$
is a holomorphic nowhere vanishing $(3,0)$-form on $X$, and for any scalar function $u$ on $M$, the $(1,1)$-form
\be\label{FYansatz}
\omega_u=\pi^*(e^u\hat\omega)+i\theta\wedge\bar\theta
\ee
is a conformally balanced metric on $X$.

\medskip
The first non-K\"ahler solution of the Hull-Strominger was actually found by Fu and Yau \cite{FY1, FY2} in 2006, using these fibrations and elliptic partial differential equations methods. Besides the solutions themselves, the Fu-Yau work also revealed the rich and deep structure of the Hull-Strominger system. Since then many other non-K\"ahler solutions have been found,
see e.g. \cite{AGF1,AGF2,FHP1,FeiY,FIUV,FIUV2,FY1,FTY,OUV}.

\medskip

In general, for a given elliptic system, there is an infinite number of flows that will admit it as stationary points, but most of them will not behave well for large time. Usually one has the freedom of selecting a well-behaved flow. But in the present case, the Anomaly flow was dictated by the requirement that the balanced condition is preserved, and there does not appear any other flow that meets this requirement. Thus it is particularly important to determine whether the Anomaly flow is a viable tool for solving the Hull-Strominger system. This is confirmed by the following theorem in the case of toric fibrations \cite{PPZ6}:

\begin{theorem}  
\label{toricfibrations-theorem}
Consider the Anomaly flow on the fibration $X\to Y$ constructed above,
with initial data $\omega(0)=\pi^*(M\hat\omega)+i\theta\wedge\bar\theta$, where $M$ is a positive constant. Fix $h$ on the bundle $E\rightarrow Y$ satisfying the Hermitian-Yang-Mills equation $\hat\omega\wedge F=0$. Then $\omega(t)$ is of the form $\pi^*(e^u\hat\omega)+i\theta\wedge\bar\theta$ and, assuming an integrability condition on the data (which is necessary), there exists $M_0>0$, so that for all $M\geq M_0$, the flow exists for all time, and converges to a metric $\omega_\infty$ with $(X, \pi^*E, \omega_\infty, \pi^*(h))$ satisfying the Hull-Strominger system.
\end{theorem}

\medskip

The reason why the Anomaly flow on Calabi-Eckmann-Goldstein-Prokushkin fibrations is simpler than the general case is because it descends to a flow on the base $Y$ of the form.
\bea
\label{af4}
\partial_t\omega
=
-{1\over 2\|\Omega\|_\omega}\left[{R\over 2}-|T|^2
-
{\alpha'\over 4}\sigma_2(iRic_\omega)
+2\alpha'{i\p\bar\p(\|\Omega\|_\o\rho)\over\o^2}-2{\mu\over\omega^2}\right]\omega
\eea
with an initial metric of the form $\omega(0)=M\hat\omega$. Here $\omega$ is now a Hermitian metric on the base $Y$, and $\sigma_2(\Phi)=\Phi\wedge\Phi\,\o^{-2}$ is the usual determinant of a real $(1,1)$-form $\Phi$, relative to the metric $\o$. $\rho$ and $\mu$ are given smooth (1, 1)-form and (2, 2)-form, respectively, both are determined explicitly by $\hat\o$, $F$ and $\theta$. 

\smallskip
Since $\omega$ stays conformal to $\hat\omega$, we can set $\omega=e^u\hat\omega$, and rewrite also this flow as a parabolic flow in $u$,
\bea
\label{FY_parabolic_u}
\p_t u = {1 \over 2} e^{-u}\left({\Delta_{\hat\o}} e^u + \alpha'  \hat{\sigma}_2(i \ddb u) - 2 \alpha'  {i \ddb (e^{-u} \rho) \over \hat{\o}^2} + 2{\mu\over\hat\o^2} \right)
\eea
where both the Laplacian $\Delta_{\hat\o}$ and the determinant $\hat\sigma_2$ are written with respect to the fixed metric $\hat\o$. This flow is not concave in the Hessian of $u$, and this makes it inaccessible to standard non-linear PDE techniques. In practice, we actually solve the problem by using the previous, Ricci-flow like formulation.

\medskip
Some key ingredients of the proof are the following estimates. We illustrate only features which are peculiar to the flow (\ref{af4}), the details can be found in the original paper \cite{PPZ6}.

\subsubsection{A priori estimate for $\|\Omega\|_\o$}

The first estimate is a $C^0$ estimate for $\|\Omega\|_\o=e^{-u}$, or equivalently an upper and a lower bound for $e^u$, but which has to be sharper than the usual $C^0$ estimates. This is because it does not suffice to show that the upper and lower bounds are independent of time, but also that they remain comparable to the size of the initial data.

\begin{lemma}
\label{C0} 
Assume that the flow (\ref{af4}) exists and is parabolic for $t\in [0,T)$ and starts from $u_0=\log M$. Then there exists $M_0$ such that, for all $M\geq M_0$, we have
\bea
{\rm sup}_{X\times [0,T)}e^u\leq C_1M,
\quad
{\rm sup}_{X\times [0,T)}e^{-u}\leq C_2M^{-1}
\eea
with constants $C_1,C_2$ depending only on $(Y, \hat\o)$, $\rho$, $\mu$, and $\alpha'$.
\end{lemma}

This estimate is proved by parabolic Moser iteration. The first step is to isolate from the flow a term comparable to the $L^2$ norm of the gradient of $e^{(k+1)u/2}$, for arbitrary $k$. For example, to establish the sup norm bound, we use
\bea
&&
{k\over 2}\int_X e^{(k+1)u}\{\hat\o +\alpha'e^{-2u}\rho\}\wedge i\p u\wedge\bar\p u
+
{\p\over\p t}{2\over k+1}\int_X e^{(k+1)u}{\hat\o^2\over 2!}
\nonumber\\
=
&&-{k\over 2}\int_X e^{ku}i\p u\wedge\bar\p u\wedge\o'
+
\int_Xe^{ku}\mu
-\alpha'(1-{1\over 1-k})\int_X e^{(k-1)u}i\p\bar\p\rho
\eea
with a form $\o'$ which is positive as long as the flow exists and is parabolic. To prove the lemma, we assume first that ${\rm sup}_X e^{-u(\cdot,t)}< \delta$ for all $t\in [0,T)$, where $\delta$ is a fixed small positive constant depending only on the data $(Y,\hat\o)$, $\rho$, $\mu$, and $\alpha'$. We later show that ${\rm sup}_X e^{-u(\cdot,t)}< \delta$ is preserved. The above estimate implies the bound
\bea
{k\over 2}\int_X e^{(k+1)u}|\na u|^2
+{\p\over\p t}{2\over k+1}\int_Xe^{(k+1)u}
\leq
(\|\mu\|_{L^\infty}+2\|\alpha'\rho\|_{C^2})
\left(\int_X e^{ku}+\int_X e^{(k-1)u}\right).
\eea
We can then apply the Sobolev inequality and iterate the bound. But the need for obtaining precise bounds in terms of $M$ requires that we need to consider not just large time, but also small time. This is because {\it some} finite bound is automatic for small time from parabolicity, but this is not the case for {\it the specific} bound of the type stated in the lemma.

\smallskip
For the lower bound, we need to use the upper bound to establish the preliminary integral estimate
\bea
\int_X e^{-u}\leq 2C_0M^{-1}
\eea
again assuming to begin with that ${\rm sup}_X e^{-u(\cdot,t)}< \delta$ for all $t\in [0,T)$. Isolating the gradient, applying the Sobolev inequality, and iterating as before starting from the above integral estimate, separately for the cases of small time and large time, we obtain
\bea
\label{C0a}
{\rm sup}_{X\times [0,T)}e^{-u}\leq C_2M^{-1}
\eea
assuming that ${\rm sup}_X e^{-u(\cdot,t)}< \delta$ for all $t\in [0,T)$.
The proof of Lemma \ref{C0} is then completed by the following argument. Start from $u(0)=\log M$ with $M$ large enough so that the condition ${\rm sup}_X e^{-u(\cdot,0)}< \delta$ is satisfied. The estimates which have already been obtained imply that this condition is preserved as long as the flow exists, and actually the better estimate (\ref{C0a}) holds. Q.E.D.

\subsubsection{A priori estimate for the torsion}

The previous $C^0$ estimate is the only one which relies on the formulation of the Anomaly flow as a scalar parabolic equation. From now on, the other estimates make use rather of the formulation (\ref{af4}) as a geometric flow of metrics.

\begin{lemma}
\label{C1}
There exists $M_0 > 1$ with the following property.
Let the flow (\ref{af3}) start from a constant function $u(0) = \log\,M$ with $M \geq M_0$. If
\bea
|\alpha'Ric_\o| \leq 10^{-6}
\eea
along the flow, then there exists $C_3 > 0$ depending only on $(Y, \hat \o)$, $\rho$, $\mu$, and $\alpha'$, such that
\bea
|T|^2 \leq C_3M^{-1}.
\eea
\end{lemma}

\smallskip
The proof is by the maximum principle applied to the function
\bea
G=\log |T|^2-\Lambda \log\|\Omega\|_\o
\eea
with $\Lambda=1+2^{-3}$. Set
\bea
F^{p\bar q}=g^{p\bar q}+\alpha'\|\Omega\|_\o^3\tilde\rho^{p\bar q}
-
{\alpha'\over 2}(Rg^{p\bar q}-R^{p\bar q})
\eea
where the $(1,1)$-form $\tilde \rho^{p\bar q}$ is defined by picking up the coefficient of $u_{\bar q p}$ in the expansion 
\bea
-\alpha'i\p\bar\p (e^{-u}\rho)=
(\alpha'e^{-u}\tilde \rho^{p\bar q}u_{\bar q p}+\cdots){\hat\o^2\over 2}.
\eea
We can then show that, at a maximum point of $G$, we must have
\bea
0\leq F^{p\bar q}\p_p\p_{\bar q}G
-
{1\over 200}|T|^2
+
C\|\Omega\|_\o(1+{\|\Omega\|_\o^{1\over 2}\over |T|}).
\eea
If $|T|\leq \|\Omega\|_\o^{1\over 2}$, the desired estimate follows is a trivial consequence of Lemma \ref{C0}. Otherwise, we obtain at a maximum point
\bea
|T|^2\leq 200 \,C\|\Omega\|_\o(1+{\|\Omega\|_\o^{1\over 2}\over |T|})
\eea
which again implies the desired estimate.

\subsubsection{A priori estimate for the curvature}

Perhaps the most technically elaborate estimate is the following $C^2$ estimate:

\begin{lemma}
\label{C2}
Start the flow with a constant function $u_0 = \log\,M$. There exists $M_0> 1$
such that for every $M\geq M_0$, if
\bea
\|\Omega\|_\o\leq C_2 M,
\quad
|T|^2 \leq C_3M^{-1}
\eea
along the flow, then
\bea
|\alpha'Ric_\o| \leq 
C_5
M^{-{1\over 2}} 
\eea
where $C_5$ only depends on $(Y, \hat \o)$, $\rho$, $\mu$, and $\alpha'$. Here, 
$C_2$ and $C_3$ are the constants given
in Lemmas \ref{C0} and \ref{C1} respectively.
\end{lemma}

This lemma is proved by considering the following quantity
\bea
|\alpha' Ric_\o|^2+\Theta |T|^2
\eea
where $\Theta$ is a suitable constant. This quantity starts from the value $0$ at time $t=0$. If we choose $\Theta$ to be
\bea
\Theta={\rm max}\{4|\alpha'|^{-1},
8(\alpha')^2(C_4+1)\}
\eea
then we can show that 
\bea
\p_t(|\alpha' Ric_\o|^2+\Theta |T|^2)\leq 0
\eea
the first time that it reaches the value $(2\Theta C_3+1)M^{-1}$. This implies that it never exceeds this value, and Lemma \ref{C2} is proved.

\subsubsection{The long-time existence of the flow}

A priori estimates for the derivatives of the curvature and of the torsion, in terms of suitable powers of $M$ and constants depending only on the geometric data can be formulated and established in the same way, when $|\alpha' Ric_\o|\leq \delta_1$ and $|T|\leq\delta_2$ for some fixed small positive constants $\delta_1$ and $\delta_2$..

\smallskip
Once the a priori estimates are available, the long-time existence of the flow can be established as follows. If we start from $u(0)=\log M$, we have $|\alpha' Ric_\o|=|T|=0$ at time $t=0$. So the set of times $t$ for which $|\alpha'Ric_\o|\leq \delta_1$,
$|T|\leq \delta_2$ on $X\times [0,t)$ is closed and not empty. By the a priori estimates, we obtain even better estimates $|\alpha'Ric_0|<C_5M^{-{1\over 2}}$, $|T|^2
\leq C_3M^{-1}$, which show that it is also open. Thus we obtain estimates for all derivatives and all time, which implies the existence of the flow for all time.

\subsubsection{The convergence of the flow}

The preceding bounds imply the existence of a sequence of times for which the flow is convergent. To get the full convergence, we introduce
\bea
J(t)=\int_X v^2\,{\hat\o^2\over 2!}
\eea
with $v=\p_t e^u$. Using the a priori estimates and the fact that $v$ has average $0$ at all times, we can show that
\bea
{dJ(t)\over dt}\leq -\eta J
\eea
for some strictly positive constant $\eta$ if $M$ is large enough. Thus $v$ tends to $0$ in $L^2$ norm, and it is not difficult, using the a priori estimates, to deduce then that it tends to $0$ in $C^\infty$ norm. The proof of Theorem \ref{toricfibrations-theorem} is complete.

\subsection{The Anomaly Flow on hyperk\"ahler fibrations over Riemann surfaces}

\smallskip

Next we describe a very recent work of T. Fei, Z. Huang, and S. Picard \cite{FHP1,FHP2}. Their set-up is certain hyperk\"ahler fibrations over Riemann surfaces originating from works of Fei \cite{Fei1,Fei2, Fei3}, which are themselves built on earlier constructions of Calabi \cite{Cal} and Gray \cite{Gray}.

\medskip

Let $\Sigma$ be a Riemann surface and a holomorphic map $\varphi:\Sigma\to {\bf CP}^1$ with $\varphi^*O(2)=K_\Sigma$. By pulling back sections of $O(2)$, we obtain three holomorphic $(1,0)$-forms $\mu_1,\mu_2,\mu_3$. Let $\hat\omega$ be the metric proportional to $i\sum_{j=1}^3\mu_j\wedge \bar\mu_j$ normalized so that $\int_\Sigma\hat\omega=1$. 

\smallskip

Next, take a flat $4$-torus $(T^4,g)$ which is viewed as a hyperk\"ahler manifold with complex structures $I,J,K$ satisfying $I^2=J^2=K^2=-1$ and $IJ=K$. Let $\omega_I,\omega_J,\omega_K$ be the corresponding K\"ahler forms. For any point $(\alpha,\beta,\gamma)\in S^2$, there is a compatible complex structure $\alpha I+\beta J+\gamma K$ on $T^4$. Let now 
$X=\Sigma\times T^4$, with the complex structure $j_\Sigma \oplus (\alpha I+\beta J+\gamma K)$, where we have let
$\varphi=(\alpha,\beta,\gamma)$.

\smallskip
Then $X$ is a non-k\"ahler complex manifold with trivial canonical bundle, and non-vanishing holomorphic $(3,0)$-form given by
$$
\Omega=\mu_1\wedge \omega_I+\mu_2\wedge\omega_J+\mu_3\wedge\omega_K
$$
Furthermore, for any real function $f\in C^\infty(\Sigma)$, as shown in \cite{Fei1,Fei2}, the $(1,1)$-form
\bea
\label{F}
\omega_f
=e^{2f}\hat\omega+e^f\omega',
\quad
{\rm with }\  \ \omega'=\alpha\omega_I+\beta\omega_J+\gamma\omega_K
\eea
is a conformally balanced metric on $X$. Furthermore, $\|\Omega\|_{\o_f}=e^{-2f}$ and 
\be \label{gcg-balanced}
\| \Omega \|_{\omega_f} \omega_f^2 = 2 {\rm vol}_{T^4} + 2 e^f \hat{\omega} \wedge \omega',
\ee

\medskip

The Anomaly flow on $X$ preserves the above ansatz and reduces then to the following conformal flow on $\Sigma$ for a metric $\omega(t)$
\bea
\label{gcg-af}
\p_t e^f=
{1\over 2}[\hat g^{z\bar z}\p_z\p_{\bar z}(e^f+{\alpha'\over 2}\kappa e^{-f})
-\kappa (e^f+{\alpha'\over 2}\kappa e^{-f})]
\eea
where $\kappa$ is the Gauss curvature of $\hat\o$. Explicitly, $\kappa=-\varphi^*\o_{FS}/\hat\o$ and is always $<0$, except at the branch points of $\varphi$, where it vanishes. Despite the extent to which Riemann surfaces have been studied, the flow (\ref{gcg-af}) does not seem to have appeared before. It may be useful for future investigations to have the following more intrinsic formulation of it,
\bea
\partial_t\omega
=
{\omega\over|\mu|^2}(-R+|\nabla\log |\mu|^2|^2)
+{\alpha'\over 2}(\Delta |\nabla\varphi|^2-|\nabla\varphi|^4)
\eea
where all norms and operators are taken with respect to the evolving metric $\omega(t)$. The following criterion for the long-time existence of the flow
in this setting was proved in \cite{FHP2}:

\medskip
\begin{theorem} 
Assume that $\alpha'>0$. Suppose that a solution to (\ref{gcg-af}) exists on $[0,T)$. If 
\be
\sup_{\Sigma \times [0,T)} \| \Omega \|_{\omega_f} = \sup_{\Sigma \times [0,T)} e^{-2f} < \infty,
\ee
then the flow can be extended to $[0,T+\epsilon)$ for some $\epsilon>0$.
\end{theorem}

The proof is short and can be given here. Let the function $u$ be defined by
\bea
u=e^f+{\alpha'\over 2}\kappa e^{-f}.
\eea
Then $u$ evolves by
\be \label{gcg-u-evol}
\p_t u = {1 \over 2} (1 - {\alpha' \over 2} \kappa e^{-f}) (\Delta_{\hat{\omega}} u - \kappa u). 
\ee
Since $\alpha'>0$, $\kappa \leq 0$, and $e^{-f} \leq C$, this is a parabolic equation for $u$ with bounded coefficients. We may therefore apply the theorem of Krylov-Safonov \cite{KrSa} to obtain parabolic H\"older estimates for $u$. Since
$2\, e^f=(u + \sqrt{u^2 - 2 \alpha' \kappa})$,
we may deduce H\"older estimates for $e^f$ and $e^{-f}$. Applying parabolic Schauder theory and a bootstrap argument, we obtain higher order estimates for $e^f$ which allow us to extend the flow past $T$. Q.E.D.

\medskip
However, it may happen that $e^{-f}\to \infty$ as $t\to T$ for some finite time $T$, and the flow will terminate in finite time. Indeed, integrating (\ref{gcg-af}) gives
\be
{d \over dt} \int_\Sigma e^f \, \hat{\omega} = {1 \over 2} \int_\Sigma (-\kappa) e^f \hat{\omega}- {\alpha' \over 4} \int_\Sigma \kappa^2 e^{-f} \hat{\omega}.
\ee
Applying the Cauchy-Schwarz inequality,
\be \label{gcg-af-ode}
{d \over dt} \int_\Sigma e^f \, \hat{\omega} \leq {\| \kappa \|_{L^\infty(\Sigma)} \over 2} \int_\Sigma e^f \hat{\omega}- {\alpha' \over 4} \left( \int_\Sigma (-\kappa) \hat{\omega} \right)^2 \left( \int_\Sigma e^f \hat{\omega} \right)^{-1}.
\ee
Since the genus $g \geq 3$ for this construction, the Gauss-Bonnet theorem gives $\int_\Sigma (-\kappa) \hat{\omega} > 0$. By studying this ODE for $\int_\Sigma e^f \hat{\omega}$, we see that if $\int_\Sigma e^f \hat{\omega}$ is too small initially then it will collapse to zero in finite time. This provides an example of the Anomaly flow on $[0,T)$ with $T<\infty$ where $\| \Omega \|_{\omega_f} = e^{-2f} \rightarrow \infty$ as $t \rightarrow T$.

\medskip
At the other end, it was proved in \cite{FHP2} that for initial data with $e^f$ large, then the flow will exist for all time, and collapse the fibers in the limit:

\begin{theorem} Assume that $e^f\gg1$ at the initial time. Then the flow exists for all time, and
$$
{\omega_f\over \int_X\|\Omega\|_{\omega_f}{\omega_f^3\over 3!}}\to
p^*(q_1^2\hat\omega)
$$
where $q_1$ is the first eigenfunction of the operator $-\Delta_{\hat\omega}-|\nabla\varphi|_{\hat\omega}^2$. Furthermore
$$
(X,{\omega_f\over \int_X\|\Omega\|_{\omega_f}{\omega_f^3\over 3!}})
\ \rightarrow\ (\Sigma, q_1^2\,\hat\omega)
\quad  {\rm in\  Gromov-Hausdorff\ topology}.
$$
\end{theorem}

\medskip
The two regions of initial data $e^f$ small or large leave out the particularly interesting region of intermediate $e^f$.
It is actually in this intermediate region that Fei, Huang, and Picard \cite{FHP1} found infinitely many topologically distinct examples, which is remarkable, as the analogous statement has been conjectured for Calabi-Yau manifolds, but has not been proved as yet.

\bigskip
When comparing the Anomaly flow to the K\"ahler-Ricci flow in section \S \ref{comparisonRF}, we had noted that that the Anomaly flow preserves the $(2,2)$ de Rham cohomology class $[\|\Omega\|_\o\o^2]$ (when $c_2(X)=c_2(E)$), but that this condition is much weaker than the preservation of the $(1,1)$-K\"ahler class as in the case of the K\"ahler-Ricci flow. Nevertheless, it can give important information, as we shall now discuss in the present case.

\smallskip
Taking the de Rham class of (\ref{gcg-balanced}) yields
\be 
{1 \over 2} [\| \Omega \|_{\omega_f} \omega_f^2] = [{\rm vol}_{T^4}] + [ e^f \alpha \hat{\omega}] [\omega_I] + [e^f  \beta  \hat{\omega}] [\omega_J] + [e^f  \gamma \hat{\omega}] [\omega_K].
\ee
Therefore $[ \| \Omega \|_{\omega_f} \omega_f^2] \in H^4(X,{\bf R})$ is parametrized by the vector 
\be
V = (\int_\Sigma e^f \alpha \, \hat{\omega}, \int_\Sigma e^f \beta \, \hat{\omega}, \int_\Sigma e^f \gamma \, \hat{\omega}).
\ee
Since the class $[ \| \Omega \|_{\omega_f} \omega_f^2]$ is preserved along the Anomaly flow, the vector $V$ is constant along the flow, and 
since $(\alpha,\beta,\gamma) \in S^2$, we have
\be
\int_\Sigma e^f \hat{\omega} \geq |V|.
\ee
This is a bound from below for the flow in terms of cohomological data.
It gives more refined information about the finite time singularity obtained when $\int_\Sigma e^f \hat{\omega} \ll 1$. Since the ODE (\ref{gcg-af-ode}) forces $\int_\Sigma e^f \hat{\omega} \rightarrow 0$ in finite time, we see that the long-time existence criterion must fail first, and $e^{-2f} = \| \Omega_f \|_{\omega_f} \rightarrow \infty$ as $t \rightarrow T$ on a set of measure zero.

\smallskip
Finally, we note that the flow admits an energy functional
\be
I(u) = \int_\Sigma |\nabla u|^2_{\hat{\omega}} \, \hat{\omega} +  \int_\Sigma \kappa u^2 \hat{\omega}.
\ee
which is monotone decreasing for arbitrary initial data. This functional was used in \cite{FHP2} to provide convergence of the normalized solution in the regime $u \geq 0$. 

It is possible that the existence of solutions to the Hull-Strominger system may be tied to some stability condition. In this case, the functional $I(u)$ may turn out to be useful. We note that \cite{FHP1}
introduced a condition which they call the ``hemisphere condition", which may also be related to some notion of stability.

\subsection{The case of unimodular Lie groups}

Next, we consider the Anomaly flow on unimodular Lie groups. This is a natural setting to consider, since any left-invariant metric will be conformally balanced. The left-invariance property will also reduce the Anomaly flow to a system of ordinary differential equations. This will allow us to tackle the first case when the Anomaly flow is a genuine system instead of just a scalar equation. It will be then easier to trace the effects of a non-zero parameter $\alpha'$. An interesting feature of the Hull-Strominger equations will emerge, which is that solutions may require a different connection on the Gauduchon line of metrics associated to a Hermitian metric than the Chern unitary connection. This had been anticipated in the physics literature \cite{AGF3, BM, FIUV, Hull1}, and worked out explicitly by Fei and Yau \cite{FeiY} for stationary points. The only caveat is that we shall be dealing with non-compact settings, which may not behave exactly in the same way as in the compact case.

\medskip

Consider left-invariant metrics on a $3$-dimensional complex Lie group $X$. Let $\{e_a\}$ be a basis of left-invariant holomorphic vector fields on $X$, with structure constants $c^d{}_{ab}$,
$$
[e_a,e_b]=\sum_d e_d c^d{}_{ab}
$$
Then $\Omega=e^1\wedge e^2\wedge e^3$ is a holomorphic nowhere vanishing  $(3,0)$-form on $X$, and if $X$ is unimodular in the sense that
$$
\sum_d c^d{}_{db}=0, \qquad {\rm for\ any}\ b,
$$
then any metric $\omega=ig_{\bar ab}e^b\wedge\bar e^a$ is balanced \cite{AG}. Unimodular Lie groups in 3-d are given by ${\bf C}^3$, the Heisenberg group, the rigid motions of ${\bf R}^2$, and $SL(2,{\bf C})$.

\medskip

Because $\omega$ is not necessarily K\"ahler, there are many natural unitary connections associated to it. The most familiar one is the Chern unitary connection, defined earlier.
But other unitary connections $\nabla^{(\kappa)}$, called the Gauduchon line, can be defined by
$$
\nabla_j^{(\kappa)}V^k=\nabla_jV^k-\kappa T^k{}_{jm}V^m
$$
Here $T=i\partial\omega$ is the torsion.
The connections with $\kappa=1$ and $\kappa=1/2$ are known respectively as the Bismut connection and the Lichnerowicz connection. The following theorem was proved in \cite{PPZ7}:

\begin{theorem} 
\label{Lie} Set $\tau=2\kappa^2(2\kappa-1)$, and assume that $\alpha'\tau>0$. Then

\smallskip

{\rm (a)} When $X={\bf C}^3$, any metric, and hence the flow is stationary.

{\rm (b)} When $X$ is nilpotent, there is no stationary point, and hence the flow can never converge. If the initial metric is diagonal, the metric remains diagonal along the flow, the lowest eigenvalue remains constant, while the other two eigenvalues tend to $\infty$ at a constant rate.

{\rm (c)} When $X$ is solvable, the stationary points of the flow are precisely the metrics with
$$
g_{\bar 12}=g_{\bar 21}=0,
\qquad \alpha'\tau\,g^{3\bar 3}=1.
$$
The flow is asymptotically unstable near any stationary point. However, the condition $g_{\bar 12}=g_{\bar 21}=0$ is preserved along the flow, and the flow with any initial metric satisfying this condition converges to a stationary point.

{\rm (d)} When $X=SL(2,{\bf C})$, there is a unique stationary point, given by the diagonal metric
$$
g_{\bar ab}={\alpha'\tau\over 2}\delta_{ab}$$
The linearization of the flow at the fixed point admits both positive and negative eigenvalues, and hence the flow is asymptotically unstable.

\end{theorem}

\medskip
The details can be found in \cite{PPZ7}. Here we shall stress only that, underlying all this is the remarkable fact, worked out by Fei and Yau, that while $Rm$ is not a $(1,1)$-form for general $\nabla^{(\kappa)}$ connections,
${\rm Tr}(Rm\wedge Rm)$ is a $(2,2)$-form. More precisely,
$$
{\rm Tr}(Rm\wedge Rm)_{\bar k\bar\ell ij}
=
\tau \overline{c^r{}_{k\ell}c^s{}_{rp}}c^q{}_{ij}c^s{}_{qp}
$$

\smallskip
We should also say that the equations appear complicated, for example, they are given by
$$
\partial_tg_{\bar 11}={1\over 2\|\Omega\|}g^{3\bar 3}g_{\bar 11}(1-{\alpha'\tau\over 4}g^{3\bar 3}),
\quad
\partial_t g_{\bar 22}={1\over 2\|\Omega\|}g^{3\bar 3}g_{\bar 22}(1-{\alpha'\tau\over 4}g^{3\bar 3}), etc.$$
in the case of the group of rigid motions, and in the case of $SL(2,{\bf C})$ by
$$
\partial_t\lambda_1
=
{(\lambda_1\lambda_2\lambda_3)^{1\over 2}\over 2}
\bigg\{
{\alpha'\tau\over 4}({2\over\lambda_1}+{\lambda_1\over\lambda_2^2}+{\lambda_1\over\lambda_3^2})-{\lambda_2\over\lambda_3}-{\lambda_3\over\lambda_2}\bigg\}
$$
with similar equations for $\lambda_2$ and $\lambda_3$. But they reveal a lot of structure upon closer inspection, and this leads to the above detailed information on the phase diagram of the corresponding dynamical systems.

\section{A Flow of Balanced Metrics with K\"ahler Fixed Points}
\setcounter{equation}{0}

The special case $\alpha'=0$ of the Anomaly flow may actually be of independent geometric interest as well. It is given by $\p_t(\|\Omega\|_\o\o^2)=i\p\bar\p\o$ for $3$-folds $X$, and a natural generalization to arbitrary $n$ dimensions for $n\geq 3$ can be defined by \cite{PPZ9}
\bea
\label{bflow}
\partial_t(\|\Omega\|_\omega\omega^{n-1})
=i\partial\bar\partial \omega^{n-2}.
\eea
Here $X$ is assumed to be a compact complex manifold of dimension $n$, which admits a nowhere vanishing holomorphic $(n,0)$-form $\Omega$, and the initial metric $\o(0)$ is assumed to be conformally balanced, $d(\|\Omega\|_{\o(0)}\o^{n-1}(0))=0$.

\subsection{Formulation as a flow of metrics}\label{balance-atheno-kahler}

As in the case of $3$ dimensions, the above flow of $(n-1,n-1)$-forms admits a
useful reformulation as a flow of metrics, as it was given in \cite{PPZ9}:
\begin{theorem}
Consider the flow (\ref{bflow}) with $X$, $\Omega$, and $\omega(0)$ as described above. Set $\o=ig_{\bar kj}dz^j\wedge d\bar z^k$. Then the evolution of $g_{\bar kj}$ is given by

{\rm (a)} When $n=3$,
\bea
\label{bflow3}
\p_tg_{\bar kj}={1\over 2\|\Omega\|_\o}\left[-\tilde R_{\bar kj}+g^{m\bar\ell}g^{s\bar r}
T_{\bar rmj}\bar T_{s\bar\ell\bar k}\right]
\eea

{\rm (b)} When $n\geq 4$,
\bea
\label{bflow4}
\p_t g_{\bar{k} j} &=& {1 \over (n-1) \| \Omega \|_\omega} \bigg[ -\tilde{R}_{\bar{k} j} + {1 \over 2 (n-2)} (|T|^2- 2 |\tau|^2) \, g_{\bar{k} j}\nonumber\\
&&- {1 \over 2} g^{q \bar{p}} g^{s \bar{r}} T_{\bar k q s} \bar{T}_{j\bar p\bar r}
+ g^{s\bar r} (T_{\bar k j s} \bar{T}_{\bar r} + T_s \bar{T}_{j\bar k \bar r})  +  T_j \bar{T}_{\bar k} \bigg].
\eea
\end{theorem}

Recall that the torsion terms are defined by
$T=i\partial\omega$,
$T_\ell=g^{j\bar k}T_{\bar kj\ell}$,
and we have introduced the $(1,0)$-form $\tau=T_\ell \,dz^\ell$.

The reason that the flow (\ref{bflow}) may be interesting is because its stationary points are Ricci-flat K\"ahler metrics. In fact, the stationary points of the flow satisfy the so-called {\it astheno-K\"ahler} condition $i\p\bar\p\o^{n-2}=0$ introduced in \cite{JY}. The fact that conformally balanced and astheno-K\"ahler metrics must be K\"ahler had been shown earlier by Matsuo-Takahashi \cite{MT} and Fino-Tomassini \cite{FT}.
With the above formulas (\ref{bflow3}) and (\ref{bflow4}) for the flow (\ref{bflow}),
it is easy to give an independent proof of this fact. We give the proof when $n\geq 4$, the same argument works for $n=3$ and is even simpler. At a stationary point $\p_tg_{\bar kj}=0$, and taking the trace of the right hand side of (\ref{bflow4}) gives
$$
0 =  -\tilde{R} + {n \over 2 (n-2)} (|T|^2- 2 |\tau|^2) - {1 \over 2} |T|^2 + 3 |\tau|^2 $$
But for conformally balanced metrics, 
$\tilde{R} = R=\Delta \log \| \Omega \|^2_\omega$. Simplifying yields
$$
(n-2) \Delta \log \| \Omega \|^2_\omega = |T|^2 + 2(n-3) |\tau|^2\geq 0$$
Thus $\log \| \Omega \|^2_\omega$ is constant by maximum principle,
$|T|^2=0$, and $\omega$ is K\"ahler and Ricci-flat.

\subsection{A convergence theorem}

It is a well-known and still open question in complex geometry to determine when a balanced or conformally balanced complex manifold $X$ is actually K\"ahler. The flow (\ref{bflow}) appears particularly well-suited to address this question when $c_1(X)=0$. In particular, from the above discussion, the K\"ahler property is a necessary condition for the convergence of the flow. At this moment, we can prove a partial converse \cite{PPZ9}:
 
\begin{theorem}
\label{bflow-theorem}
Assume that the initial data $\omega(0)$ satisfies 
$\| \Omega \|_{\omega(0)} \omega(0)^{n-1} = \hat{\chi}^{n-1}$,
where $\hat{\chi}$ is a K\"ahler metric. Then the flow exists for all time $t>0$, and as $t \rightarrow \infty$, the solution $\omega(t)$ converges smoothly to a K\"ahler, Ricci-flat, metric $\omega_\infty$ satisfying
$$
\omega_\infty = \| \Omega \|_{\chi_\infty}^{-2/(n-2)} \chi_\infty$$
where $\chi_\infty$ is the unique K\"ahler Ricci-flat metric in the cohomology class $[\hat{\chi}]$, and $\| \Omega \|_{\chi_\infty}$ is an explicit constant.
\end{theorem}

\bigskip

It is instructive to examine the corresponding PDE in greater detail.

\medskip

The key starting point is that, in this special case, we can actually identify an evolution of $(1,1)$-cohomology classes. Define the smooth function 
$f\in C^\infty(X,{\bf R})$ by
$$
e^{-f} = {1 \over (n-1)}  \| \Omega \|^{-2}_{\hat{\chi}}$$
Let $t\to\varphi(t)$ be the following
Monge-Amp\`ere type flow
$$
\p_t \varphi = e^{-f} {(\hat\chi+i\partial\bar\partial\varphi)^n
\over\hat\chi^n}, \quad {\rm with } \ \ \varphi(x,0)=0,$$
subject to the plurisubharmonicity condition $\hat{\chi} + i \ddb \varphi >0$.
Then, $\o(t)$ defined by 
\be\label{ansatz}
\| \Omega \|_{\omega(t)}\, \omega^{n-1}(t) = \chi^{n-1}(t) \quad {\rm with } \ \ \chi(t) = \hat{\chi} + i \ddb \varphi(t) > 0
\ee
is precisely the solution of the Anomaly flow.

\medskip

We note the similarity with the K\"ahler-Ricci flow and with the flow recently introduced by Collins et al \cite{Co1}, except we have here the Monge-Amp\`ere determinant itself, instead of its log or its inverse.
In this respect, the Anomaly flow with $\alpha'=0$ can be viewed more like the complex analogue of the inverse Gauss curvature flow, which for convex bodies can be expressed as
$$
\p_t u
= {\det (u g_{ij}+ \nabla_i \nabla_j u ) \over \det g_{ij}}, \ \ \ u(x,0)=u_0(x)>0$$
where $u$ is the support function $u: {\bf S}^n \times [0,T)$ defined by $u(N,t) = \langle P, N \rangle$ where $P \in M_t$ is the point on $M_t$ with normal $N$, and $({\bf S}^n,g_{ij})$ is the standard sphere.

\medskip

But unlike the K\"ahler-Ricci flow, the equation is not concave. Thus many well-known methods for the theory of non-linear parabolic partial differential equations, including those recently developed in \cite{PT}. Rather, the equation has to be studied on its own. We describe now some of the main steps.

\subsubsection{Estimates for $\p_t\varphi$}

Despite the difficulties we have noted, the flow has a very important property, which is that $\p_t\varphi\equiv \dot\varphi$ as well as the oscillation of $\varphi$ can be controlled at all times. To see this, let $L$ be the linear elliptic second order operator defined by
\bea
L=e^{-f}{(\hat \chi+i\p\bar\p\varphi)^n\over\hat\chi^n}\chi^{j\bar k}\hat\na_j\hat\na_{\bar k}
\eea
where we have denoted by $\hat\na$ the connection with respect to the reference metric $\hat\chi$. Then differentiating the flow gives
\bea
(\p_t-L)\dot\varphi=0
\eea
and hence, by the maximum principle
\bea
{\rm inf}_X\dot\varphi(0)
\leq \dot\varphi\leq {\rm sup}_X\dot\varphi(0).
\eea
which shows that $\dot\varphi$ is uniformly bounded for all times. Next, we can rewrite the flow as a Monge-Amp\`ere equation with bounded right hand side,
\bea
(\hat\chi+i\p\bar\p\varphi)^n=e^f\dot\varphi \hat\chi^n
\eea
Yau's $C^0$ estimate gives then a uniform bound on the oscillation
${\rm sup}_X\varphi-{\rm inf}_X\varphi$.

\subsubsection{The $C^2$ estimate}

The key estimate is the following $C^2$ estimate: let $\varphi$ be a solution of (\ref{bflow}) on $X\times [0,T]$ satisfying the positivity condition (\ref{ansatz}). Then there exists constants $C, A>0$ depending only on $X$, $\hat\chi$, and $f$ such that
\bea
\Delta_{\hat\chi}
\varphi(x,t)\leq C\,e^{A(\tilde \varphi(x,t)-{\rm inf}_{X\times [0,T]}\tilde\varphi)}
\eea
for all $(x,t)\in X\times [0,T]$. Here $\tilde\varphi$ is defined by
\bea
\tilde\varphi(x,t)=\varphi(x,t)-{1\over V}\int_X \varphi(\cdot,t)\hat\chi^n,
\qquad
V=\int_X\hat\chi^n.
\eea

The proof is by the maximum principle applied to the following test function
\bea
G(x, t) = \log {\rm Tr} \,h - A \tilde \varphi +
{B\over 2}[{(\hat\chi+i\p\bar\p\varphi)^n\over\hat\chi^n}]^2
\eea
where $h^j{}_k=\hat\chi^{j\bar p}\chi_{\bar pk}$, and $A,B > 0$ are suitable constants. We note the extra
term $B[(\hat\chi+i\p\bar\p\varphi)^n/\hat\chi^n]^2/2$, which does not appear in the standard test function used in the study of the
complex Monge-Amp\`ere type equations or the K\"ahler-Ricci flow. Indeed, this is the main innovation of this test function, and it is required in order to overcome
the difficult terms caused by the lack of concavity of the flow as well as cross terms resulting from the conformal factor $e^{-f}$.

\subsubsection{The $C^{2,\alpha}$ estimate}

The lack of concavity is again a source of difficulties in getting $C^{2,\alpha}$ estimates and the absence of an Evans-Krylov theorem. Here we adapt instead techniques for $C^{2,\alpha}$ estimates developed by B. Chow and D. Tsai \cite{CT, Tsai} (see also Andrews \cite{Andrews}) in the study of Gauss curvature type flows in convex geometry. This implies the long-time existence of the flow.

\smallskip
Finally the convergence of the flow can be established in several ways. One elegant way is to use the dilaton functional
\bea
\label{M}
M(\o)=\int_X \|\Omega\|_\o\o^n
\eea
introduced by Garcia-Fernandez et al \cite{GFRST}. This functional can be shown to be monotone decreasing for the flow (\ref{bflow}) with conformally K\"ahler initial data, with $dM(t)/dt\to 0$ as $t\to\infty$. This completes the proof of Theorem \ref{bflow-theorem}.

\subsection{Singularities and divergence}

Since the convergence of the flow (\ref{bflow}) is intimately related to the question of whether the underlying conformally balanced manifold $X$ is actually K\"ahler, it would be very valuable to find geometric criteria for when it develops singularities and/or diverges. An elementary example is provided by the case of Calabi-Eckmann-Goldstein-Prokushkin fibrations with $\alpha'=0$. If $\o_1$ and $\o_2$ are the two harmonic forms used in the construction of the toric fibration $\pi: X\to Y$, and $\o_u$ is the metric (\ref{FYansatz}), then 
$$i\p\bar\p \o_u=(\Delta_{\hat\o}e^u+\|\o_1\|_{\hat\o}^2+\|\o_2\|_{\hat\o}^2)\hat\o^2/2,$$ so that the Anomaly flow
$\p_t(\|\Omega\|_{\o_u}\o_u^2)=i\p\bar\p\o_u$ becomes
\bea
\p_t e^u
=
{1\over 2}(\Delta_{\hat\o}e^u+\|\o_1\|_{\hat\o}^2+\|\o_2\|_{\hat\o}^2).
\eea
It follows immediately that
\bea
\p_t\left(\int_Y e^u\hat\o^2\right)
=\int_Y \left(\|\o_1\|_{\hat\o}^2+\|\o_2\|_{\hat\o}^2\right){\hat\o^2\over 2}>0.
\eea
Thus the flow cannot converge unless $\o_1=\o_2=0$, which is exactly when the fibration $\pi$ becomes the trivial product $Y\times T$, which is indeed K\"ahler.

\medskip
The study of the formation of singularities in a flow begins with the identification of which geometric quantity must blow up in order for the flow not to extend further. For the original Anomaly flow in dimension $3$, the following theorem was established in \cite{PPZ5}: 

\begin{theorem} Assume that $n=3$ and we have a solution of the Anomaly flow (\ref{af3}) with $\alpha'=0$ (or (\ref{bflow3})) on the time interval $[0,{1\over A}]$. Then for all non-negative integers $k$, there exists a constant $C_k$, depending on a uniform lower bound for $\|\Omega\|_\o$, such that, if
\bea
|Rm|_\o+|DT|_\o
+|T|_\o^2\leq A, \qquad (z,t)\in X\times[0,{1\over A}]
\eea
then
\bea
|D^kRm(z,t)|_\o+|D^{k+1}T(z,t)|_\o\leq {C_kA\over t^{k\over 2}},
\qquad
(z,t)\in X\times (0,{1\over A}].
\eea
In particular, the flow will exist for all time, unless there is a time $T>0$ and a sequence $(z_j,t_j)$ with $t_j\to T$ with either $\|\Omega(z_j,t_j)\|_{\o(t_j)}\to 0$, or
\bea
\left(|Rm|_\o^2+|DT|_\o^2+|T|_{\o}^2\right)(z_j,t_j)\to\infty.
\eea
\end{theorem}

Indeed, the above theorem holds for the flow (\ref{bflow3}) in any dimension $n\geq 3$. We expect this theorem to generalize to the flow (\ref{bflow}) for arbitrary $n\geq 3$.

\section{More Speculative Flows}
\setcounter{equation}{0}

We conclude with a brief discussion of some more speculative flows. 

\subsection{Flow of balanced metrics with source term}

In \S 6, we investigated the flow (\ref{bflow}) for conformally balanced metrics, which is just the higher dimensional generalization of the Anomaly flow (\ref{af1}) with zero slope parameter $\alpha'$. We can also consider a similar flow with source term
\be\label{flow-source}
\p_t(\|\Omega\|_{\o} \omega^2) = i\ddb \omega - \Upsilon,
\ee
where $\Upsilon$ is a given $i\p\bar\p$-exact (2,2)-form on the compact complex $3$-fold $X$. Indeed, this flow provides a natural approach to a question raised by Garcia-Fernandez \cite{GF}, namely, whether for any given real $i\p\bar\p$-exact $(2,2)$-form $\Upsilon$ on a compact conformally balanced $3$-fold $(X,\Omega,\o_0)$, there exists another conformally balanced metric $\o$ with
$i\p\bar\p\o= \Upsilon$ and $[\|\Omega\|_\o\o^2]=
[\|\Omega\|_{\o_0}\o_0^2]$ ? This question it itself closely related 
to a conjecture on the solvability of the Hull-Strominger system, raised by Yau. Thus the long time behavior and convergence of the flow (\ref{flow-source}) is of particular interest.

\subsection{The Anomaly flow in arbitrary dimension}

The original motivation for the Hull-Strominger system, and hence for the Anomaly flow was for compactifications of the $10$-dimensional heterotic string to $4$-dimensional space-times, and thus the internal complex manifold $X$ must have complex dimension $3$. The Green-Schwarz anomaly cancellation mechanism made also use specifically of the second Chern classes of $T^{1,0}(X)$ and of the bundle $E\to X$, which appear in (\ref{anomaly}). However, from the pure mathematical standpoint, the Anomaly flow seems to fit naturally in a family of flows which can be defined in arbitrary dimension $n$, and are given by
\bea
\partial_t(\|\Omega\|_\omega\omega^{n-1})
=
i\partial\bar\partial\omega^{n-2}-\alpha'(c_{n-1}(Rm)-\Phi_{n-1}(t))
\eea
where $c_{n-1}(Rm)={\rm Tr} (Rm\wedge\cdots\wedge Rm)$ is the $(n-1)$-th Chern class of $T^{1,0}(X)$, and $\Phi_{n-1}(t)$ is a given flow of closed $(n-1,n-1)$-forms, for example $\Phi_{n-1}(t)=c_{n-1}(F(t))$ where $F(t)$ is the curvature of a flow of Hermitian metrics $H(t)$.

\smallskip

Just as Hermitian-Yang-Mills equations originally appeared in dimension $4$ but proved to be mathematically important in arbitrary dimensions, we expect that these flows will be interesting for $n>3$ as well. We note that other generalizations of the Hull-Strominger system to arbitrary dimensions have been proposed earlier in \cite{FY1, FY2}. Restricted to Calabi-Eckmann-Golstein-Prokushkin  fibrations, they lead to the following interesting $n$-dimensional generalizations of the equation studied by Fu and Yau \cite{FY1, FY2}, called the Fu-Yau equation,
\be \label{FY-form-eq}
i \ddb (e^u \hat{\omega} - \alpha' e^{-u} \rho) \wedge \hat{\omega}^{n-2} + \alpha' i \ddb u \wedge i \ddb u \wedge \hat{\omega}^{n-2} + \mu = 0,
\ee
where $\mu$ is a given smooth $(n, n)$ form. When $n=2$, the above equation reduces to the right hand side of equation (\ref{FY_parabolic_u}) and was studied by Fu and Yau using elliptic methods in \cite{FY1, FY2}. The general dimensional case of (\ref{FY-form-eq}) has been subsequently studied in \cite{PPZ1, PPZ4, PPZ8, CHZ1, CHZ2}. They are of Hessian type, and suggest interesting questions about a priori estimates for such equations, some of which have been addressed in \cite{DK, Guanbo, GRW, PPZ3, PPZ4}.

\subsection{Flows without the form $\Omega$}

For some questions, e.g., the relation between balanced cone and K\"ahler cone \cite{FX}, it would be desirable to eliminate the assumption of the existence of a non-vanishing holomorphic form $\Omega$. Thus it is tempting to consider the flow
$$
\partial_t\omega^{n-1}
=
i\partial\bar\partial \omega^{n-2}$$
on an $n$-dimensional complex manifold $X$, for initial data $\o(0)$ which are balanced, i.e. $d\o(0)^{n-1}=0$. The stationary points are still exactly the K\"ahler metrics, although they are not necessarily Ricci-flat. Explicitly, we can still write it in terms of metrics,
\bea
\p_tg_{\bar kj}
=-{1\over n-1}\na^mT_{\bar kjm}+
{1\over 2(n-1)}\bigg[-g^{q\bar p}g^{s\bar r}T_{\bar kqs}\bar T_{j\bar p\bar r}
+{|T|^2\over n-1}g_{\bar kj}\bigg]
\eea
Compared to the expression (\ref{bflow4}) for the Anomaly flow, we see that the key term $\tilde R_{\bar kj}$ has cancelled out, and the flow is not parabolic. The lack of parabolicity poses severe problems, and it is unclear whether this flow can be used in any way. It would be very interesting if one can find a flow of balanced metrics with K\"ahler metrics as fixed points, but which is parabolic.
In this context, we note that another flow of balanced metrics has been introduced by Bedulli and Vezzoni in \cite{BV}, based on a generalization of the Calabi flow, and its short-time existence established there.

\

\bigskip

Department of Mathematics, Columbia University, New York, NY 10027, USA

\smallskip

phong@math.columbia.edu

\smallskip
Department of Mathematics, Columbia University, New York, NY 10027, USA

\smallskip
 picard@math.columbia.edu

\smallskip
Department of Mathematics, University of California, Irvine, CA 92697, USA

\smallskip
xiangwen@math.uci.edu

\end{document}